\documentclass{amsart}
\usepackage{pstricks}
\usepackage{epsfig}
\usepackage{pst-grad} 
\usepackage{pst-plot} 

\usepackage{amsmath}%
\usepackage{amsfonts}%
\usepackage{amssymb}%
\usepackage{graphicx}

\newtheorem{theorem}{Theorem}
\theoremstyle{plain}

\newtheorem{corollary}{Corollary}

\newtheorem{definition}{Definition}

\newtheorem{lemma}{Lemma}

\newtheorem{proposition}{Proposition}

\numberwithin{equation}{section}
\newcommand{\R}{\mathbb{R}}
\newcommand{\C}{\mathbb{C}}
\DeclareMathOperator{\imag}{Im}
\DeclareMathOperator{\ind}{Ind}

\DeclareMathOperator{\spn}{span}
\begin{document}
\title[Global $\widetilde{SL(2,R)}$ representations of the Schr\"{o}dinger equation]{Global $\widetilde{SL(2,R)}$ representations of the Schr\"{o}dinger equation with time-dependent potentials}
\author[Jose Franco]{Jose Franco\\
Department of Mathematics \\ Baylor University \\ One Bear Place \# 97329 \\ Waco, TX\\ jose\_franco@baylor.edu}%
\subjclass{22E70, 35Q41} %
\keywords{Schr\"{o}dinger equation, time-dependent potentials, Lie theory, representation theory, globalizations}%

\begin{abstract}
We study the representation theory of the solution space of the one-dimensional Schr\"{o}dinger equation with time-dependent potentials that posses $\mathfrak{sl}_2$-symmetry. We give explicit local intertwining maps to multiplier representations and show that the study of the solution space for potentials of the form $V(t,x)=g_2(t)x^2+g_1(t)x+g_0(t)$ reduces to the study of the potential free case. We also show that the study of the time-dependent potentials of the form $V(t,x)=\lambda x^{-2}+g_2(t)x^2+g_0(t)$ reduces to the study of the potential $V(t,x)=\lambda x^{-2}$. Therefore, we study the representation theory associated to solutions of the Schr\"{o}dinger equation with this potential. The subspace of solutions for which the action globalizes is constructed via nonstandard induction outside the semisimple category.
\end{abstract}
\maketitle

\section{Introduction}

In the early seventies, the original prolongation algorithm of Sophus Lie was used to classify the various time-independent potentials of the 
one-dimensional Schr\"{o}dinger equation, $$2iu_t+u_{xx}=2V(t,x)u$$  that admit non-trivial, inequivalent Lie symmetries (c.f. \cite{Boyer}, \cite{Niederer}). As local Lie group representations, it turns out that the representation theory associated to the space of solutions of the Schr\"{o}dinger equation with the following potentials, 
\begin{subequations}\label{Potentials}
\begin{gather}
V_1(x)=\lambda\label{ConstPot}\\
V_2(x)=\lambda x\label{LinearPot}\\
V_3(x)=\lambda x^2\label{QuadPot},
\end{gather}
\end{subequations}
where $\lambda\in \R$ is an arbitrary constant, reduces to the potential free case. Here, the symmetry algebra is isomorphic to $\mathfrak g :=\mathfrak{sl}(2,\R)\ltimes \mathfrak{h}_3(\R)$. More generally, it is known that the time-dependent potential $V(t,x)=g_2(t)x^2+g_1(t)x+g_0(t)$ has the same symmetry Lie algebra $\mathfrak{g}$ (c.f. \cite{Truax1}).

However, for the following potentials,
\begin{subequations}\label{Potentials2}
\begin{gather}
V_4(x)=\lambda x^{-2}\label{InvQuadPot}\\
V_5(x)=\lambda_1 x^2+\lambda_2 x^{-2} \label{ComboPot},
\end{gather}
\end{subequations}
where $\lambda, \lambda_i \in \R$ are arbitrary constants, we will show that the problem will reduce to a study of an eigenvalue problem for, essentially, a Casimir element for $\mathfrak{sl}(2,\R)$. For these cases, the symmetry Lie algebra is isomorphic to $\mathfrak{sl}(2,\R)\times \R$. Similarly, Truax showed that the time-dependent potential $V(t,x)=\lambda x^{-2}+g_2(t)x^2+g_0(t)$ has the same symmetry Lie algebra.

It is natural to use representation theory to study the solution space of these differential operators. However, since the resulting actions are not global, the techniques of representation theory do not always apply. 
However, it is sometimes possible to look at special subspaces of solution that carry the structure of a global representation (c.f. \cite{Craddock}, \cite{Sepanski1}, \cite{Sepanski2}). For instance, in 2005, M. Sepanski and R. Stanke decomposed the solution
space for the $1$-dimensional potential free Schr\"{o}dinger equation and studied it as a global Lie group representation in \cite{Sepanski2}. Recently, they analyzed the $n$-dimensional case for the potential free Schr\"{o}dinger equation (c.f. \cite{Sepanski1}). 

In this paper, we give explicit local intertwining maps to multiplier representations, showing that the study of the solution space for the potentials of the form $V(t,x)=\lambda x^{-2}+g_2(t)x^2+g_0(t)$, which includes \eqref{ComboPot}, reduces to the study of the potential \eqref{InvQuadPot}. For the sake of completeness, we show that the study of potentials  \eqref{Potentials} and the time-dependent potentials of the form $V(t,x)=g_2(t)x^2+g_1(t)x+g_0(t)$, reduces to the study of the potential free case. Therefore, we study the representation theory associated to solutions of the Schr\"{o}dinger equation with potentials \eqref{InvQuadPot}. As in \cite{Sepanski1} and \cite{Sepanski2}, the subspace of solutions for which the action globalizes is constructed via nonstandard induction outside the semisimple category. 

A bit more precisely, we start with a parabolic-like subgroup $\overline{P}$ of the group $G:=\widetilde{SL(2,\R)}\ltimes H_3$ where $\widetilde{SL(2,\R)}$ is the two-fold cover of $SL(2,\R)$ and $H_3$ is the three-dimensional Heisenberg group. Then, we look at the smoothly induced representation $$I(q,r,s)=\ind_{\overline P}^{G}(\chi_{q,r,s})$$ where $\chi_{q,r,s}$ is a character on $\overline{P}$ with parameters $r,s \in \C$ and $q \in \mathbb Z_4$ (see Equation \eqref{Character}). 

We show that in the non-compact version of $I(q,r,s)$ for $r=-1/2$ and $s=i/2$, solutions to the Schr\"{o}dinger equation with potential \eqref{InvQuadPot} are realized as solutions to the eigenvalue problem \begin{equation}\label{eigenProb}\Omega f=\left(\frac{r(r+2)}{2}+\lambda\right)f\end{equation}
where $\Omega$ is the Casimir element of $\mathfrak{sl}(2,\R)$. We show that the solution space in $I(q,r,s)$ is non-empty iff $\lambda = l(l-1)/2$ for some $l \in \mathbb Z^{\geq 0}$. 

To each $l \in \mathbb Z^{\geq 0}$ we assign the triangular number $\lambda = l(l-1)/2$. This relation is one-to-one except for $l =0$ and $l=1$. Under this identification, we denote by $H_l$ the space of $K$-finite vectors of the solution space in $I(q,r,s)$ of \eqref{eigenProb} for $l\geq 2$. For $\lambda = 0$ we decompose the space of $K$-finite vectors of the solution space into two $\mathfrak{sl}_2$-invariant subspaces $H_0\oplus H_1$ (see Theorem \ref{BigTheorem}).

We determine the structure of $H_l$ as an $\mathfrak{sl}(2,\R)$-module and show:

\begin{theorem}
Assume $\lambda = l(l-1)/2$ for some $l \in \mathbb Z^{\geq 0}$.  If $q = 0$ or $q=2$  then $H_l$ is irreducible as an $\mathfrak{sl}_2$-module. If $q= 1$ (respectively, $q =3$), then $H_l$ has a unique lowest (respectively, highest) weight submodule, $H_l^-$ (respectively, $H_l^+$).
\end{theorem}

While the action of the Heisenberg group does not preserve $H_l$ we show that it does preserve the direct sum $H:=\bigoplus_{l=0}^\infty H_l$. Also, write $H^\pm$ for the direct sum of highest and lowest weight submodules with $l\geq 2$, whenever they exist.  We prove:

\begin{theorem} As $\mathfrak{g}$-modules:
\begin{enumerate}
	\item If $q = 0$ or $q = 2$, the composition series of $H$ is $$0\subset H_0\oplus H_1 \subset H.$$
	\item If $q = 1$, then the composition series of $H$ is as follows  $$0\subset H_0^-\oplus H_1^- \subset H_0\oplus H_1\subset H_0\oplus H_1\oplus H^-\subset H.$$
	\item If $q =3$, then then the composition series of $H$ is  $$0\subset H_0^+\oplus H_1^+ \subset H_0\oplus H_1\subset H_0\oplus H_1\oplus H^+\subset H.$$
\end{enumerate}
\end{theorem}

\section{Notation}

\subsection{The Group}

Let $G_0=SL(2,\R)$ and let $H_3$ denote the three dimensional Heisenberg group with product,
$$(v_1,v_2,v_3)(v'_1,v'_2,v'_3) = (v_1+v'_1,v_2+v'_2,v_3+v'_3 + v_1v'_2 - v_2v'_1). $$
Following the realization of the two-fold cover of $G_0$ in \cite{Kashi}, define the complex upper half plane $D:=\{z\in\C|\imag z >0\}$ and let $G_0$ act on $D$ by fractional linear transformations, that is, if $g=\left(\begin{smallmatrix} a & b \\ c & d \end{smallmatrix}\right)\in G_0$ and $z\in D$, then $$g.z=\frac{az+b}{cz+d}.$$ Define $d:G_0\times D \to \C$ by $d(g,z):=cz+d$. Then there are exactly two smooth square roots of $d(g,z)$ for each $g \in G_0$ and $z\in D$. The double cover can be realized as:
\begin{multline*}
\widetilde{G_0}=\left\{(g,\epsilon) | g\in SL(2,\R) \text{ and smooth }\epsilon:D\to\C \right. \\ \left. \text{ such that } \epsilon(z)^2=d(g,z)\text{ for } z\in D \right\}
\end{multline*}
with the product defined by $$(g_1,\epsilon_1(z))(g_2,\epsilon_2(z))=(g_1g_2, \epsilon_1(g_2.z)\epsilon_2(z)).$$
Finally, the symmetry group that we are interested in, is $G:=\widetilde{G_0}\ltimes H_3$. Here $\widetilde G_0$ projects to $G_0$ and acts on $H_3$ by the standard action on the first two coordinates and leaves the third fixed.
\subsection{Parabolic Subgroup and Induced Representations}
As in \cite{Sepanski1}, we consider the parabolic subalgebra of lower triangular matrices $\overline{\mathfrak{q}} \subset \mathfrak{sl}(2,\R)$ with Langlands decomposition $\mathfrak{m \oplus a} \oplus \overline{\mathfrak{n}}$. If $\exp_{\widetilde G_0}: \mathfrak{sl}(2,\R)\to \widetilde{G_0}$ denotes the exponential map then:
\begin{gather*}
A:= \exp_{\widetilde{G_0}}(\mathfrak a)=\{(\left(\begin{smallmatrix}t & 0 \\ 0&t^{-1}\end{smallmatrix}\right),z\mapsto e^{-t/2})| t\in\R^{\geq 0}\}\\
N:= \exp_{\widetilde{G_0}}(\mathfrak n)=\{(\left(\begin{smallmatrix}1 & t \\ 0& 1\end{smallmatrix}\right), z\mapsto 1)|t\in\R\}\\
\overline{N}:= \exp_{\widetilde{G_0}}(\overline{\mathfrak{n}})=\{(\left(\begin{smallmatrix}1 & 0 \\ t & 1\end{smallmatrix}\right), z\mapsto \sqrt{tz+1})|t\in\R\}.
\end{gather*}
Let $\mathfrak{k}:= \{\left(\begin{smallmatrix}0 & \theta \\ -\theta & 0 \end{smallmatrix}\right):\theta \in \R\}$ then $$K:= \exp_{\widetilde{G_0}}(\mathfrak k)=\{ (\left(\begin{smallmatrix}\cos\theta & \sin\theta \\ -\sin\theta & \cos\theta\end{smallmatrix}\right),z \mapsto \sqrt{\cos\theta-z\sin\theta} )|\theta \in \R\}$$
where $\sqrt{\cdot}$ denotes the principal square root in $\C$. Writing $M$ for the centralizer of $A$ in $K$ then
$$M=\{m_j:=(\left(\begin{smallmatrix}-1 & 0 \\0 &-1\end{smallmatrix}\right)^j, z\to i^{-j}) | j=0,1,2,3\}.$$
Let $W \subset H_3$ be given by $W=\{(0,v,w)| v,w \in \R \}\cong \R^2$ and let $X:=\{(x,0,0) |x \in \R\}$. Let us write $\mathfrak{w}$ for the Lie algebra of $W$. Then $\overline P=MA\overline N \ltimes W$ is the analogue of a parabolic subgroup in $G$ corresponding to $\overline{\mathfrak p}:=\overline{\mathfrak q} \ltimes \mathfrak w$. 

For later use, we notice that an element in $g=\left[(\left(\begin{smallmatrix} a& b \\ c & d\end{smallmatrix}\right),z\mapsto \epsilon(z)),\left( 
u,v,w\right) \right]\in G$ is in the image of the mapping $\overline{P}\times(N\times X)\to G$ given by $(\overline{p},n)\mapsto \overline{p}n$, if $a\neq 0$. This induces a decomposition of such $g$ into its $\overline P$ and $N\times X$ components,
\begin{multline*}\left[(\left(\begin{smallmatrix} a& b \\ c & d\end{smallmatrix}\right),z\mapsto \epsilon(z)),\left( 
u,v,w\right) \right]=\\ \left[(\left(\begin{smallmatrix} a& 0 \\ c & a^{-1}\end{smallmatrix}\right),z\mapsto \epsilon(z+b/a)),\left( 
0, v, w + (u + bv/a)v\right) \right] \\ \cdot \left[(\left(\begin{smallmatrix} 1& b/a \\ 0 & 1\end{smallmatrix}\right),z\mapsto 1),\left( 
u+bv/a,0,0\right) \right].
\end{multline*}
On the open dense set where $a\neq 0$ let $\overline{p}:G\to \overline{P}$ and $n:G\to N\times X$  be the projections from the previous decomposition.

It is well known that the character group on $A$ is isomorphic to the additive group $\C$ so any character on $A$ can be indexed by a constant $r \in \C$ and defined by 
$$\chi_r\bigl((\left(\begin{smallmatrix}t & 0 \\ 0&t^{-1}\end{smallmatrix}\right),z\mapsto e^{-t/2})\bigr)=t^r$$
for $t>0$. A character on $M$ is parametrized by $q\in \mathbb Z_4$ and defined by $\chi_q(m_j)=i^{j q}$. A character on $W$ can be parametrized by $s\in \C$ and defined by,
$$\chi_{s}\bigl((0,v,w)\bigr)=e^{sw}.$$
 Finally, any character on $\overline P$ that is trivial on $N$ is parametrized by a triplet $(q,r,s)$ where $s, r \in \C$ and $q\in \mathbb Z_4$ and defined by
\begin{equation}\label{Character}\chi_{q,r,s}\left(((-1)^j\left(\begin{smallmatrix} a& 0 \\ c & a^{-1}\end{smallmatrix}\right),z\mapsto i^{-j}e^{-a/2}\sqrt{acz+1}),\left( 
0,v,w\right) \right)=i^{jq}|a|^r e^{sw}.\end{equation}
The representation space induced by $\chi_{q,r,s}$ will be denoted by $I(q,r,s)$ and defined by 
$$I(q,r,s):= \{\phi:G\to\C | \phi \in C^\infty \text{ and } \phi(g\overline p)=\chi_{q,r,s}^{-1}(\overline p)\phi(g) \text{ for } g\in G, \overline p\in \overline P\}\\$$
the $G$-action on $I(q,r,s)$ is given by $(g_1.\phi)(g_2)=\phi(g_1^{-1}g_2)$.

\subsection{Non-compact Picture}
Since $H_3=XW$ then $G=(N\times X)\overline P$ a.e. with $N\times X$ isomorphic to $\R^2$ via $(t,x)\mapsto N_{t,x}:=\left[(\left(\begin{smallmatrix}1 & t \\ 0& 1\end{smallmatrix}\right), z\mapsto 1),(x,0,0)\right]$. Since a section in the induced representation is determined by its restriction to $N$, this restriction induces an injection of $I(q,r,s)$ into $C^\infty(\R^2)$ which is identified as
\begin{equation*}I'(q,r,s)=\{f\in C^\infty(\R^2) | f(t,x)=\phi(N_{t,x})  \text{ for some } \phi\in I(q,r,s)\}.\end{equation*} 

This space is endowed with the corresponding action so that the map $\phi \mapsto f$ where $f(t,x)=\phi(N_{t,x})$, becomes intertwining. Thus $I(q,r,s)\cong I'(q,r,s)$ as $G$-modules. As in the semisimple case, we will call this the non-compact picture. 

\begin{proposition}\label{GroupActionStandard}
Let $f \in I'(q,r,s)$, $(g,\epsilon)\in \widetilde{G_0}$, and $(u,v,w)\in H_3$. Then,
\begin{subequations}
 \begin{align}
((g,\epsilon).f)(t,x)&=(a-ct)^{r-q/2} \epsilon(g^{-1}.(t+z)) e^{\frac{-s c x^2}{a-ct}} f\left(\frac{dt-b}{a-ct},\frac{x}{a-ct} \right) \label{SL2StdAction}\\
 ((u,v,w).f)(t,x)&=e^{-s(uv-2vx-tv^2+w)}f(t,x-u-tv). \label{H3StdAction}
\end{align}
\end{subequations}
\begin{proof}
This result is proved in a more general setting in \cite{Sepanski1}.
\end{proof}
\end{proposition}
\begin{corollary}\label{sl2ActionsStd}
The action of $\left(\begin{smallmatrix}a & b \\ c & -a\end{smallmatrix}\right)\in\mathfrak{sl}(2,\R)$ on $I'(q,r,s)$ is given by the differential operator
\begin{equation}\label{AlgActionNonComp}(ct-a)x\partial_x+(ct^2-2at-b)\partial_t+(ra-csx^2-rct). \end{equation}
An element $(u,v,w)\in \mathfrak{h}_3$ acts on $I'(q,r,s)$ by the differential operator
$$(tv-u)\partial_x+s(w-2vx).$$
\begin{proof}It follows from differentiating the group actions on $I'(q,r,s)$. 
\end{proof}
\end{corollary}

\subsection{Casimir Elements}

Write $$\Box = 2i\partial_t+\partial_x^2$$ for the potential free Schr\"{o}dinger operator. By equation \eqref{AlgActionNonComp}, the standard $\mathfrak{sl}_2$-triple $\{h,e^\pm\}$ acts by
\begin{gather}
h= -x\partial_x-2t\partial_t+r\\
e^+= -\partial_t\\
e^-=tx\partial_x+t^2\partial_t-(sx^2+rt)
\end{gather}
on the non-compact picture. Let  $$\Omega=1/2h^2-h+2e^+e^-$$ be the Casimir element in the enveloping algebra of $\mathfrak{sl}(2,\R)$ and define 
$$\Omega'=2\Omega-r(r+2).$$

\begin{corollary}\label{CasimirNonComp}
On $I'(q,r,s)$, $\Omega$ acts by $$ \Omega = \frac{1}{2}\left(4sx^2\partial_t+x^2\partial_x^2-(1+2r)x\partial_x+r(r+2) \right).$$
In particular, for $r=-1/2$ and $s=i/2$, $\Omega'$ acts by
$$\Omega'- 2\lambda = x^2 (\Box-2\lambda/x^2)$$
so that,
$$\ker(\Omega'- 2\lambda)=\ker(\Box-2\lambda/x^2).$$
\begin{proof}A straightforward calculation using the actions of the standard $\mathfrak{sl}_2$-triple and the definition of the Casimir element gives the desired result.
\end{proof}
\end{corollary}
As a consequence of Corollary \ref{CasimirNonComp}, we are interested in the study of $\ker(\Omega'- 2\lambda)$. We begin with a result on the invariance under $\widetilde{G_0}$ and the subgroup $\{(0,0,w)|w\in\R\}\subset H_3$.

\begin{proposition}
The subspace $\ker(\Omega'-2\lambda)$ in $I'(q,r,s)$ is invariant under the action of $\widetilde{G_0}$ and under the action of the subgroup $\{(0,0,w)|w\in\R\}$ of $H_3$. The space is not left invariant by the complement of $\{(0,0,w)|w\in\R\}$ in $H_3$.
\begin{proof}
Since $\Omega$ is in the center of the enveloping algebra of $\mathfrak g$, the $\widetilde{G_0}$-invariance is clear. Let $(u,v,w)\in H_3$. Using the action on Corollary \ref{sl2ActionsStd} we can calculate
\begin{equation*}
\begin{split}[\Box-2\lambda/x^2,(u,v,w)]&=[2i\partial_t+\partial_x^2-2\lambda/x^2,(tv-u)\partial_x+i/2(w-2vx)]\\&  =[2i\partial_t,(tv-u)\partial_x]+[\partial_x^2,i/2(w-2vx)]\\ & \quad -2[\lambda/x^2,(tv-u)\partial_x]\\ 
& = 2iv\partial_x-2iv\partial_x+\frac{4\lambda(tv-u)}{x^3}=\frac{4\lambda(tv-u)}{x^3}.
\end{split}
\end{equation*}
\end{proof}
\end{proposition}

Though, all of $H_3$ does not leave $\ker(\Omega'-2\lambda)$ invariant, it will play an important role in linking together different $\widetilde G_0$-invariant kernels. 
\section{Multiplier Representations}

In this section, we consider the potential $V(t,x)=g_2(t)x^2+g_1(t)x+g_0(t)$. We construct local intertwining maps from the non-compact picture, $I'(q,r,s)$ to a multiplier representation. This will allow us to reduce the study to the potential free case.

It was shown in \cite{Truax1} that the only time-dependent potentials having full $\mathfrak{sl_2}$-symmetry are potentials of the form $V(t,x)=g_2(t)x^2+g_1(t)x+g_0(t)+\lambda/x^2$ with $\lambda \cdot g_1(t)=0$. If $\lambda=0$, the symmetry Lie algebra is isomorphic to $\mathfrak g=\mathfrak{sl}_2\ltimes \mathfrak{h}_3$ and it is isomorphic to $\mathfrak{sl}(2,\R)\times \R$ otherwise.

When $\lambda=0$ we show that there exists a local intertwining isomorphism between the solution space of the potential free Schr\"{o}dinger equation and the solution space of the Schr\"{o}dinger equation with this time-dependent potential in $I'(q,r,s)$.  Notice that the cases where $g_i(t)=\lambda$ and $g_j(t)=0$ for $i\neq j$ are the time-independent potentials listed in \eqref{Potentials}. 

When $\lambda \neq 0$ we show that, in the same multiplier representation, the study reduces to the eigenvalue problem that is the main concern of this paper. Therefore, this completes the study of all non-trivial time-dependent and independent potentials having at least the $\mathfrak{sl}_2$-symmetry.

To write down the algebra generators explicitly, fix two real, nontrivial, linearly independent solutions of $b''+2g_2b=0$. These solutions, $\chi_1$ and $\chi_2$, are normalized so that $W(\chi_1,\chi_2)=1$. With $\chi_1$ and $\chi_2$, we define $\varphi_j=\chi_j^2$ for $j \in \{1,2\}$ and $\varphi_3=2\chi_1\chi_2$. Let $\mathcal C_j(t)=\int_0^t\chi_j g_1$ for $j \in \{1,2\}$, $$\mathcal A_l = -\chi_l\mathcal C_l$$ for $l \in \{1,2\}$, and $$\mathcal A_3 = -(\chi_1\mathcal C_2+\chi_2\mathcal C_1).$$ 
It was shown in \cite{Truax1} that the differential operators:
\begin{equation}\label{sl2generators}L_j=(-1)^{j+1}(\varphi_j\partial_t+\bigl(\frac{1}{2}\varphi_j'x+\mathcal A_j\bigr)\partial_x +\mathcal B_j)\end{equation} 
for $1\leq j\leq 3$, generate an algebra isomorphic to $\mathfrak{sl}(2,\R)$, where $$\mathcal B_j = -\frac{1}{4}i\varphi_j'' x^2 -i \mathcal A'_j x+\frac{1}{4}\varphi'_j+ig_0\varphi_j + i\mathcal D_j, $$
$\mathcal D_l = -\frac{1}{2}\mathcal C_l^2$ for $l \in \{1,2\}$, and $\mathcal D_3= -\mathcal C_1\mathcal C_2$. The following bracket relations hold: $[L_3,L_1]=-2L_1$, $[L_3,L_2]=2L_2$, and $[L_2,L_1]=L_3$. Note that, in order to have a standard $\mathfrak{sl}_2$-triple, our definition of $L_2$ differs in sign with respect to the definition in \cite{Truax1}.

In order to define the appropriate multiplier representation space we start by defining a change of variables $\gamma : \R^2 \to \R^2$ by 
$$\gamma(t,x):=\left(\int_0^t\frac{1}{\chi_2^2} \ ,\ \frac{1}{\chi_2(t)}x+\int_0^t\frac{\mathcal C_2}{\chi_2^2} \right).$$
In the following, we assume that the required integrability conditions are satisfied.

We define $\nu: N\times X\to \C$ by $$\nu(N_{\gamma(t,x)})= e^{\int_0^t \frac{\mathcal B_2(u)}{\chi^2_2(u)}+\bigl(\frac{1}{\chi^2_2(u)}(\frac{1}{2}\varphi_2'(u)x+\mathcal A_2(u))\bigr)^2du}$$ and extend it to a map on an open dense subset of $G$ by $\nu(g)=\nu(n(g)).$ Let $f\in I'(q,r,s)$ and define the map $f\mapsto \tilde{f}$ by 
\begin{equation}\label{TimeDepChange}\tilde{f}(t,x)=e^{\int_0^t \frac{\mathcal B_2(u)}{\chi^2_2(u)}+\bigl(\frac{1}{\chi^2_2(u)}(\frac{1}{2}\varphi_2'(u)x+\mathcal A_2(u))\bigr)^2du}f(\gamma(t,x)).\end{equation}
The space $I'(q,r,s)_\mu$ is defined as the image of $I'(q,r,s)$ under this map (the reason for the subscript $\mu$ will become evident below). This space is provided with the structure of a $G$-module that makes the map intertwining. 

Next we construct the multiplier representation. We start by defining the multiplier
$$\mu(g_1,g_2)=\nu(g_2^{-1}g_1)\nu(g_2^{-1})^{-1}.$$ For $\phi \in I(q,r,s)$ define, on an open dense set of $G/\overline P$, the map $$\tilde \phi(g\overline P)=\mu(g^{-1},I)^{-1}\phi(g).$$ We define $I(q,r,s)_\mu$ as the image of $I(q,r,s)$ under the map $\phi \to \tilde \phi$.

Finally, the intertwining map from $I(q,r,s)_\mu$ to $I'(q,r,s)_\mu$ is given by 
\begin{equation*}\tilde\phi \mapsto \tilde f \text{ whenever }\tilde{f}(t,x)=\tilde\phi(N_{\gamma(t,x)}\overline P).\end{equation*}
\subsection{Group action on $I'(q,r,s)_\mu$} In this section we calculate the local action of $G$ on $I'(q,r,s)_\mu$ and we show that the study of the solution space for this general potential reduces to the study of the kernel of the differential operator $\Omega'$ as in the potential free case. For notational convenience, define 
$$\rho(t,x)=\nu(N_{\gamma(t,x)}).$$
\begin{proposition}
Fix $g=\left(\begin{smallmatrix} a& b \\ c & d\end{smallmatrix}\right)\in G_0$ and let $(g,\epsilon)\in \widetilde{G_0}$. Define $\Theta(t)=\int_0^t\frac{1}{\chi_2^2}$ and $\Xi(t,x)=\frac{1}{\chi_2(t)}x+\int_0^t\frac{\mathcal C_2}{\chi_2^2} $. Then
\begin{multline}\label{SL2ActionTimeDep}
((g,\epsilon).\tilde{f})(t,x) =  \frac{\rho(t,x)}{\rho\circ\gamma^{-1}\left(\frac{d\Theta-b}{a-c\Theta},\frac{\Xi}{a-c\Theta} \right)}(a-c\Theta)^{r-q/2}\\ \epsilon(g^{-1}.(\Theta+z)) e^{\frac{-s c \Xi^2}{a-c\Theta}}\tilde{f}\circ \gamma^{-1}\left(\frac{d\Theta-b}{a-c\Theta},\frac{\Xi}{a-c\Theta} \right).
\end{multline}
For $(u,v,w)\in H_3(\R)$ 
\begin{multline}
((u,v,w).\tilde{f})(t,x)= \frac{\rho(t,x)}{\rho\circ\gamma^{-1}\left(\Theta,\Xi-u-v\Theta \right)}e^{-s(uv-2v\Xi-v^2\Theta+w)}\\ \tilde{f}\circ \gamma^{-1}\left(\Theta,\Xi-u-v\Theta \right).
\end{multline}
\begin{proof}
This is a straightforward calculation. It follows directly from using the isomorphism determined by Equation \eqref{TimeDepChange} on the actions computed in Proposition \ref{GroupActionStandard}.
\end{proof}
\end{proposition}
Differentiating these actions, we recover the generators of the algebra of symmetry operators found by Truax in \cite{Truax1}. We start with some useful calculations
\begin{lemma}\label{RelationsOfChis}
The functions $\chi_1(t)$ and $\chi_2(t)$ satisfy
\begin{enumerate}
	\item $$\chi_1=\chi_2 \int_0^t\frac{1}{\chi^2_2},$$
	\item $$\mathcal A _1 =\left(\int_0^t\frac{1}{\chi^2_2}\right)^2 \mathcal A _2-\chi_2\left(\int_0^t\frac{1}{\chi^2_2}\right)\int_0^t\frac{\mathcal C_2}{\chi^2_2},$$
	\item $$\varphi'_1 = \left(\int_0^t\frac{1}{\chi^2_2}\right)^2 \varphi'_2+2\int_0^t\frac{1}{\chi^2_2}.$$
\end{enumerate}
\begin{proof}
This lemma follows from the definitions together with the fact that $\chi_1\chi_2'-\chi_2\chi_1'=1$ and the chain rule.
\end{proof}
\end{lemma}
 
\begin{corollary}
For $r=-1/2$ and $s=i/2$ the standard $\mathfrak{sl}_2$-basis $\{h,e^+,e^-\}$ acts on $I'(q,r,s)_{\mu}$ by the differential operators $\{L_3,L_2,L_1\}$ respectively. An element $(u,v,w)\in H_3(\R)$ acts on the same space, by
$$(u\chi_2-v\chi_1)\partial_x-i(u\chi_2'+v\chi_1')x+i(u\mathcal C_2+v\mathcal C_1-s w)$$ 
\begin{proof}
All these calculations are similar. We only provide the details for the action of $e^-$ since it is the most involved. Let $\gamma^{-1}(t,x)=(\Theta^{-1}(t),\Psi(t,x))$. Using \eqref{SL2ActionTimeDep}, we obtain
\begin{multline*}
((\left(\begin{smallmatrix} 1& 0 \\ c & 1\end{smallmatrix}\right),\epsilon).\tilde{f})(t,x)=\frac{\rho(t,x)}{\rho\left(\Theta^{-1}\left(\frac{\Theta}{1-c\Theta}\right),\Psi\left(\frac{\Xi}{1-c\Theta}\right) \right)}(1-c\Theta)^{r-q/2}\\ \epsilon(g^{-1}.(\Theta+z)) e^{\frac{-s c \Xi^2}{1-c\Theta}}\tilde{f}\left(\Theta^{-1}\left(\frac{\Theta}{1-c\Theta}\right),\Psi\left(\frac{\Xi}{1-c\Theta}\right) \right).
\end{multline*}
We next take $\left.\frac{d}{dc}\right|_{c=0}$ to obtain the action of $e^-$. For the coefficient of $\partial_t$ we obtain
$$-\Theta^2(t)\left.\frac{\partial \Theta^{-1}}{dt}\right|_{t=\Theta}=\chi_2(t)^2\left(\int_0^t\frac{1}{\chi_2^2}\right)^2=\chi_1^2(t)=\varphi_1(t).$$
The first equality above follows from differentiating, the second equality from Lemma \ref{RelationsOfChis}, and the third equality from the definition of $\varphi_1$. 

For the coefficient of $\partial_x$ we obtain
\begin{multline*}
\Theta(t)\Xi(t,x) \frac{\partial \Psi}{\partial x}+\Theta(t)^2\frac{\partial \Psi}{\partial t} =\Theta(t)\Xi(t,x)\chi_2+\Theta(t)^2(\frac{1}{2}\varphi'_2x+\mathcal A_2)=\\ x\Theta+\frac{1}{2}x(\varphi_1'-2\Theta)+\mathcal A_1  =\frac{1}{2}\varphi'_1x+\mathcal A_1.
\end{multline*}
The first equality above follows from differentiating, the second from computing the partial derivatives of the inverse function, and the third by Lemma \ref{RelationsOfChis}. 

Lastly, we compute the multiplication term. To this end, we start by computing the following expression
\begin{multline*}
\mathcal B _2 \Theta^2-\frac{i}{2}\Xi^2-r\Theta= \mathcal B_1+\frac{ix^2\varphi_2'\Theta}{\chi_2^2}+ix(\Theta \frac{\mathcal C_2}{\chi_2}-\chi_2'\Theta \int_0^t\frac{\mathcal C_2}{\chi_2^2})-i\Theta \mathcal C_2 \int_0^t\frac{\mathcal C_2}{\chi_2^2}.
\end{multline*}
Since the multiplication term is given by $\Theta^2 \mathcal B_2-\frac{i}{2}\Xi^2-r\Theta+\frac{1}{\rho(t,x)}\frac{\partial\rho}{\partial x}\frac{\partial\Psi}{\partial x}\Xi(t,x)\Theta=\Theta^2 \mathcal B_2+i\Theta(-\chi_2'x+\mathcal C_2)(1/\chi_2 x+\int \mathcal C_2/\chi_2^2)=\mathcal B_1$, the operator corresponding to $e^-$ is $L_1$. 
\end{proof}
\end{corollary}

Using this action of the $\mathfrak{sl}_2$-triple, a straightforward but long calculation gives the following corollary.

\begin{corollary}
For the parameters $r=-1/2$ and $s=i/2$, the Casimir element acts on $I(q,r,s)_\mu$ by $$\Omega=\frac{1}{2}\left[(x-\chi_1\mathcal C_2+\chi_2\mathcal C_1)^2(\Box-2(g_2(t)x^2+g_1(t)x+g_0(t)))-3/4\right].$$
In particular, $$\ker \Omega' = \ker\left(\Box-2(g_2(t)x^2+g_1(t)x+g_0(t))\right)$$ in $I'(q,r,s)_\mu$. 

If $\lambda\neq 0$ then $g_1(t)\equiv 0$ and $\Omega'$ acts by $\Omega'=x^2(\Box-2(g_2(t)x^2+g_0(t)))$. Thus $$\ker (\Omega'-2\lambda) = \ker\left(\Box-2(g_2(t)x^2+g_0(t)+\lambda/x^2\right)$$ in $I'(q,r,s)_\mu$.

\end{corollary}

This shows that when $\lambda=0$, at least locally, the time-dependent cases analyzed reduce to the potential free case. It also shows that when $\lambda \neq 0$ the study reduces to the study of the eigenvalue problem of $\Omega'$. (See Corollary \ref{CasimirNonComp}.) 

\section{The Compact Picture and $K$-types}

The group $\widetilde{G_0}$ has Iwasawa decomposition $\widetilde{G_0}=KA\overline{N}$ and the product induces a diffeomorphism $G\cong (K\times X) \times (A\overline{N}\ltimes W)$. Since $(A\overline{N}\ltimes W)\subset \overline P$, an element $\phi \in I(q,r,s)$ is completely determined by its restriction to $K\times X$. Moreover, since $(K \times X)\cap \overline P=M$ we have that the restriction of $\phi \in I(q,r,s)$ (which we will still denote by $\phi$) satisfies $\phi(gm)=\chi_{q,r,s}(m)^{-1}\phi(g)$ for $g\in K\times X$ and $m\in M$. 

There exists an isomorphism $K\times X \cong S^1\times \R$ given by $(\theta,y)\mapsto [(g_\theta,\epsilon_\theta),(y,0,0)]$ and it can be shown that this map is $4\pi$-periodic with respect to $\theta$. Thus we can identify $\phi\in I(q,r,s)$ with a map $F:S^1\times \R\to \C$, $\phi \mapsto F$ iff $\phi([(g_\theta,\epsilon_\theta),(y,0,0)])=F(\theta,y)$.  Then $F\in C^\infty(S^1\times \R)$ and $F(\theta+4\pi,y)=F(\theta,y)$. 

The function $F$ inherits from $\phi$ additional ``parity" identities. By the definition, $\epsilon_{\theta+\pi j}(z)^2=\cos(\theta+\pi j)-z\sin(\theta+\pi j)=(-1)^j\epsilon_\theta(z)$. We then get
\begin{multline*}F(\theta+\pi j,(-1)^jy)=\phi([(g_{\theta+\pi j},\epsilon_{\theta+\pi j}),((-1)^j y,0,0)])=\phi([(g_\theta,\epsilon_\theta),(y,0,0)]m_j)\\ =\chi_{q,r,s}(m_j)^{-1}\phi([(g_\theta,\epsilon_\theta),(y,0,0)])=i^{-jq}F(\theta,y).
\end{multline*}
Define \begin{equation}\label{defI''}I''(q,r,s)=\{F\in C^\infty(\R^2) | F(\theta+j\pi,(-1)^jy)=i^{-jq}F(\theta,y) \}. \end{equation} Then
the map $\phi \mapsto F$ is a vector space isomorphism between $I(q,r,s)$ and $I''(q,r,s)$. The space $I''(q,r,s)$ inherits a unique $G$-module structure, so that this map  becomes an intertwining map. We call this the compact picture, as in the semisimple case, though $K\times X$ is not compact here.

In turn, the isomorphism $T$ induces an isomorphism between $I'(q,r,s)$ and $I''(q,r,s)$ which we will write out explicitly. We begin with the following decomposition:
\begin{multline*}[(g_\theta,\epsilon_\theta),(y,0,0)]=[(\left(\begin{smallmatrix}1 & \tan\theta \\ 0 & 1\end{smallmatrix}\right),z\mapsto 1),(y\sec\theta,0,0)]\\ \cdot [(\left(\begin{smallmatrix} 1/\cos\theta & 0 \\-\sin\theta & \cos\theta \end{smallmatrix}\right),\epsilon_\theta ),(0,-y\tan\theta,y^2\tan\theta)].\end{multline*}
Since $F(\theta,y)=\phi([(g_\theta,\epsilon_\theta),(y,0,0)])$ then  
\begin{multline}\label{CompactPicIsomorphism}F(\theta,y)=\chi_{q,r,s}([(\left(\begin{smallmatrix} 1/\cos\theta & 0 \\-\sin\theta & \cos\theta \end{smallmatrix}\right),\epsilon_\theta ),(0,-y\tan\theta,y^2\tan\theta)])^{-1} \\ \cdot \phi([(\left(\begin{smallmatrix}1 & \tan\theta \\ 0 & 1\end{smallmatrix}\right),z\mapsto 1),(y\sec\theta,0,0)])  =(\cos\theta)^{-r}e^{-sy^2\tan\theta}f(\tan\theta,y\sec\theta)
\end{multline}
for $f\in I(q,r,s)$ and $\theta \in (-\pi/2,\pi/2)$. Since $F\in I''(q,r,s)$, this expression can be extended smoothly to any $\theta\in\R$ by using the fact that $F(\theta+j\pi,y)=i^{-jq}F(\theta,(-1)^j y)$ and continuity to get to the integer multiples of $\pi/2$. Then we define the isomorphism $T:I'(q,r,s)\to I''(q,r,s)$ by $T(f)=F$.
The inverse to this map can be calculated and it is:
\begin{equation}
f(t,x)=(1+t^2)^{-r/2}e^{\frac{stx^2}{1+t^2}}F(\arctan t ,x(1+t^2)^{-1/2}).
\end{equation}
Under this isomorphism, via the chain rule, we obtain
\begin{subequations}\label{DerivativesInCompPic}
\begin{gather}
\partial_t\leftrightarrow\frac{1}{2}(-y\sin2\theta \partial_y+\cos^2\theta\partial_\theta+2sy^2\cos2\theta-1/2r\sin2\theta) \label{PartialTCompact}\\
\partial_x\leftrightarrow 2sy\sin\theta+\cos\theta \partial_y. \label{PartialXCompact}
\end{gather}
\end{subequations}
This will enable us to transfer the actions of the algebra from the non-compact picture, $I'(q,r,s)$, to the compact picture, $I''(q,r,s)$. 

Define a standard basis of $\mathfrak{sl}_2(\C)$ given by $$\kappa = i(e^--e^+)$$ and $$\eta^\pm=1/2(h\pm i(e^++e^-)).$$ Applying equations \eqref{DerivativesInCompPic} to the action in Corollary \ref{sl2ActionsStd}, it can be shown that the $\mathfrak{sl}_2$-triple just defined acts on $I''(q,r,s)$ by the differential operators
\begin{gather}
\kappa = i\partial_\theta \label{ActionKappa}\\
\eta^\pm  = \frac{1}{2}e^{\mp2i\theta}\left(y\partial_y\mp i\partial_\theta-(1/2\pm 2isy^2)\right). \label{ActionEtas}
\end{gather}

\begin{proposition}
If $\Omega ''$ denotes the differential operator by which the central element $\Omega'$ acts on $I''(q,r,s)$ then
$$\Omega''=y^2\left(4s\partial_\theta+4s^2y^2+ \partial_y^2+\frac{1+2r}{y} \partial_y\right). $$
\begin{proof}Under the isomorphism $I'(q,r,s)\cong I''(q,r,s)$ we obtain the following expressions:
$$4sx^2\partial_t \mapsto y^2(-4sy\tan\theta \partial_y+4s\partial_\theta+8s^2y^2-4s^2y^2\sec^2\theta+2sr\tan\theta) $$
$$x^2\partial_x^2\mapsto y^2(4s^2y^2\tan^2\theta+ \partial_y^2+4sy\tan\theta\partial_y+2s\tan\theta)$$
$$(1+2r)x\partial_x\mapsto(1+2r)y^2(2s\tan\theta+1/y \partial_y). $$
Adding them we get
$$\frac{1}{y^2}\Omega''=4s\partial_\theta+4s^2y^2+ \partial_y^2+\frac{1+2r}{y} \partial_y.$$
\end{proof}
\end{proposition}

\begin{lemma}\label{quickcalculation}
Let $(g_{\theta'},\epsilon_{\theta'})\in K$ and $F\in I''(q,r,s)$ then $(g_{\theta'},\epsilon_{\theta'}).F(\theta,y)=F(\theta-\theta',y)$
\begin{proof}
\begin{multline*}(g_{\theta'},\epsilon_{\theta'}).F(\theta,y)=\phi([(g_{\theta'},\epsilon_{\theta'})^{-1}(g_\theta,\epsilon_\theta),(y,0,0)])\\=\phi([(g_{\theta-\theta'},\epsilon_{\theta-\theta'}),(y,0,0)])=F(\theta-\theta',y)\end{multline*}
\end{proof}
\end{lemma}
There exists an isomorphism $K\cong S^1$ given by $(g_\theta,\epsilon_\theta)\mapsto e^{i\theta/2}$. Therefore, the characters on $K$ are of the form $\chi_m^K(g_\theta,\epsilon_\theta)=e^{-im\theta/2}$ for $m\in \mathbb Z$.  Using Lemma \ref{quickcalculation}, a weight vector $F_m\in I''(q,r,s)$ of weight $m$, for the action of $K$, satisfies $(g_{\theta'},\epsilon_{\theta'}).F_m(\theta,y)=F_m(\theta-\theta',y)=e^{-im\theta'/2}F_m(\theta,y)$. Setting $\theta = 0$ and $\theta'=-\theta$ we obtain $F_m(\theta,y)=e^{-im\theta/2}F_m(0,y)$. Let $\tilde{F}_m(y):=F_m(0,y)$ so that a weight vector is of the form $$F_m(\theta,y)=e^{-im\theta/2}\tilde{F}_m(y).$$
\begin{lemma}\label{CondOnKtypes}
Fix $m\in \mathbb Z$ and $\tilde{F}\in C^\infty(\R)$. Then $F(\theta,y)=e^{-im\theta/2}\tilde{F}(y)$ is annihilated by $\Omega''-2\lambda$ if and only if $\tilde{F}(y)$ is annihilated by the differential operator $$\mathcal D = y^2 \partial_y^2-(2\lambda-my^2+y^4)$$
\begin{proof}
Explicitly calculating the action of $\Omega''-2\lambda$ on $F(\theta,y)=e^{-im\theta/2}\tilde{F}(y)$, one obtains that $(\Omega''-2\lambda)F=e^{-im\theta/2}\mathcal D\tilde{F}.$
\end{proof}
\end{lemma}
\begin{proposition}\label{KTypesForm}
There exist a $K$-finite vector of weight $m$ in $\ker(\Omega''-2\lambda)\subset I''(q,r,s)$ iff  
\begin{equation}\label{FormOfK}l=\frac{1}{2}(1+\sqrt{1+8\lambda})\end{equation}
is a positive integer (equivalently, $\lambda = l(l-1)/2$ for $l \in \mathbb Z^{>0}$) and $m \equiv 2l+q \mod 4$. In this case, if $\lambda \neq 0$, there exists a unique (up to scalar multiples) weight vector of weight $m$ in $\ker(\Omega''-2\lambda)\subset I''(q,r,s)$ given by
\begin{equation}\label{KFinForm}F_m(\theta,y)=e^{-im\theta/2}e^{-y^2/2}y^l\phantom{.}_1F_1\left(\frac{1+2l-m}{4},l+\frac{1}{2},y^2 \right)\end{equation}
\begin{proof}
By Lemma \ref{CondOnKtypes}, for $F_m$ to be in $\ker(\Omega''-2\lambda)\subset I''(q,r,s)$, it is necessary that $\mathcal D\tilde F_m=0$. Because of the form of $\mathcal D$, it respects the decomposition of $\tilde F_m$ in terms of its even and odd components. Moreover, each of the components is determined by its value on $\R^+$. Working first with $y\geq 0$ we can write $\tilde{F}_m(y)=e^{-y^2/2}H(y^2)$ for some smooth function $H$. Then, the condition $\mathcal D \tilde{F}_m=0$ is equivalent to
\begin{equation}\label{EqnH}
\bigl(4y^4\partial_y^2+(2y^2-4y^4)\partial_y+((m-1)y^2-2\lambda))H(y^2)=0.
\end{equation}
Following \cite{Coddington}, the Frobenius method for this equation yields a solution spanned by two linearly independent solutions. The indicial roots for this equation are $$l_1=\frac{1}{2}(1-\sqrt{1+8\lambda})$$ and $$l_2=\frac{1}{2}(1+\sqrt{1+8\lambda}).$$ Then, the first linearly independent solution is of the form $$H_1(y^2)=y^{l_2}(1+\sum_{j=1}^\infty c_j(l_2) y^{2j})$$ for some $c_j(l_2)\in \R$. This function extends to a smooth function on $\R$ only if $l_2 \in \mathbb Z^{\geq 0}$ iff $\lambda=0$ or $\lambda$ is a triangular number (i.e., $\lambda = k(k-1)/2$ for some $k \in \mathbb Z^{>1}$). 

If $\lambda\neq 0$, the difference between the indicial roots is an odd integer (i.e., $\sqrt{1+8\lambda}$), and the second solution is of the form \begin{equation}\label{ExpansionSecondSolution}
H_2(y^2)=a H_1(y^2)\ln|y^2|+y^{l_1}(1+\sum_{j=1}^\infty c_j(l_1) y^{2j})
\end{equation}

for some $a, c_j(l_1)\in \R$. Since $l_1<0$, $H_2$ is not continuous at $y=0$. 

Let $l=l_2$ and write $\tilde F_m(y)=e^{-y^2/2}y^l L(y^2)$. Applying the differential operator $\mathcal D$ to a function of the form $e^{-y^2/2}y^l L(y^2)$ we obtain the differential equation \begin{equation}\label{DEHyp}4y^2L''(y^2)+2(1 + 2l - 2y^2)L'(y^2)-(1 +2 l - m)L(y^2)=0.\end{equation}
Recall that the confluent hypergeometric differential equation is $$(z\partial_z^2+(b-z)\partial_z - a)_1F_1(a,b,z)=0$$ (c.f. \cite{Abramowitz}). This equation has well known solutions in the form of confluent hypergeometric functions of the first and second kind. However, the smoothness condition required by being in $I''(q,r,s)$ shows that the unique solution to \eqref{DEHyp} corresponds to a multiple of the confluent hypergeometric function of the first kind. We may therefore take $L(y^2)=\phantom{.}_1F_1\left(\frac{1+2l-m}{4},l+\frac{1}{2},y^2 \right)$.

Finally, a simple calculation using the required parity condition on elements in $I''(q,r,s)$ from Equation \eqref{defI''} applied to $F_m(\theta,y)$ reduces to
$$e^{-im\pi j/2}(-1)^{jl}=i^{-jq}$$
which is equivalent to $m-2l\equiv q \mod 4$.

So far, we have established the theorem for non-negative values of $y$. Extend $\tilde{F}$ to $\R$ by $\tilde F_m(y)=e^{-y^2/2}y^l \phantom{.}_1F_1\left(\frac{1+2l-m}{4},l+\frac{1}{2},y^2 \right)$ which is even or odd depending on the parity of $l$. Since $\mathcal D\tilde F_m(y)=0$ for $y\geq 0$ and $\mathcal D$ is even, $\mathcal D\tilde F_m(y)=0$ for $y\in \R$. Moreover, $F_m$ is manifestly smooth and the unique extension to all $\R$.

If $\lambda=0$ then $l=0$ or $l=1$, which corresponds to the potential free case and again, it is known that there exists a unique solution for each $l$. The solutions correspond to the even ($l=0$) and the odd ($l=1$) solutions found there. (c.f., \cite{Sepanski1})
\end{proof}
\end{proposition}

Notice that we have set up a correspondence between the set of eigenvalues with non-empty eigenspace in $I''(q,r,s)$ and $\mathbb Z^{\geq 0}$ via $\lambda =l (l-1)/2$. This correspondence will be one-to-one (except for $\lambda=0$ where it is two-to-one). For $\lambda \neq 0$ the corresponding parameter $l$ can be recovered by $l =\frac{1}{2}(1+\sqrt{1+8\lambda})$. For the potential free case, $\lambda =0$, we have associated the parameters $l=0$ and $l=1$.

For use in the following section, we record the following properties of the congruent hypergeometric function (c.f. \cite{Abramowitz})
\begin{subequations}
\begin{gather}
\frac{d^n}{dz^n}\ _1F_1(a,b,z)=\frac{(a)_n}{(b)_n}\ _1F_1(a+n,b+n,z) \label{U0} \\
b\ _1F_1(a,b,z)-b\ _1F_1(a-1,b,z)-z\ _1F_1(a,b+1,z)=0 \label{U1} \\
\begin{split}
b ( 1-b+z)\,_1F_1(a,b,z)+b(b-1)\,_1F_1(a&-1,b-1,z) \\  & -az\,_1F_1(a+1,b+1,z)=0 \label{U2} \end{split} \\
\begin{split}
(a-1+z)\,_1F_1(a,b,z)+(b-a)\,_1F_1(a&-1,b,z) \\  & (1-b)\,_1F_1(a,b-1,z)=0 \label{U3} \end{split} \\
(a-b+1)\,_1F_1(a,b,z)-a\,_1F_1(a+1,b,z)+(b-1) \,_1F_1(a,b-1,z)=0 \label{U4}
\end{gather}
\end{subequations}

\section{Structure of $\ker(\Omega''-2\lambda)\subset I''(q,r,s)$}

In this section we will study the structure of $\ker(\Omega''-2\lambda)_K$ as an $\mathfrak{sl}_2$-module.

\begin{proposition}\label{ActionsOfSL2onKTypes}
With $l=\frac{1}{2}(1+\sqrt{1+8\lambda})$ as in Proposition \ref{KTypesForm} and $m\equiv 2l+ q \mod 4$, let $$\Psi_{m,l}(\theta,y)=e^{-im\theta/2}e^{-y^2/2}y^l\,_1F_1\left(\frac{1+2l-m}{4},l+\frac{1}{2},y^2 \right).$$ The $\mathfrak{sl_2}$-triple $\{\kappa, \eta^\pm\}$ acts on $\Psi_{m,l}$ by
\begin{gather} 
\kappa.\Psi_{m,l}=\frac{m}{2}\Psi_{m,l} \\
\eta^\pm.\Psi_{m,l}=-\frac{2l+1\pm m}{4}\Psi_{m\pm4,l} \label{etasAction}
\end{gather}
Lowest weight vectors occur if $m\equiv 2l+1\mod 4$ and the lowest weight vector is of the form $$ e^{-\frac{1}{2}(2l+1)i\theta}e^{-\frac{y^2}{2}}y^l.$$
Highest weight vectors occur if $m\equiv -2l-1\mod 4$ and the highest weight vector is of the form $$ e^{\frac{1}{2}(2l+1)i\theta}e^{\frac{y^2}{2}}y^l.$$
\begin{proof}
In Equations \eqref{ActionKappa} and \eqref{ActionEtas}, we wrote down the action of the $\mathfrak{sl_2}$-triple $\{\kappa, \eta^\pm\}$. The stated action of $\kappa$ follows by inspection. 
Directly applying the differential operator $\eta^+$ gives
\begin{multline*}\eta^+.\Psi_{m,l}(\theta,y)= e^{-i(m\pm4)\theta/2}e^{-y^2/2}y^l\frac{-1 - 2 l + m}{4(1+2l)}((1 + 2 l)\\ \cdot \,_1F_1\Bigl(\frac{1+2l-m}{4},l+\frac{1}{2},y^2\Bigr)+2y^2\,_1F_1\Bigl(\frac{5+2l-m}{4},l+\frac{3}{2},y^2\Bigr)).\end{multline*}
Applying \eqref{U2} with $a=\frac{1+2l-m}{4}$ and $b =l+\frac{1}{2}$ to the action of $\eta^+$ we obtain
\begin{multline*}\eta^+.\Psi_{m,l}(\theta,y)= -\frac{1}{4}e^{-i(m\pm4)\theta/2}e^{-y^2/2}y^l((4 l-2)\\ \cdot \,_1F_1\Bigl(\frac{-3+2l-m}{4},l-\frac{1}{2},y^2\Bigr)-(3-2l+m)\,_1F_1\Bigl(\frac{1+2l-m}{4},l+\frac{1}{2},y^2\Bigr)).\end{multline*}
Using \eqref{U1} we obtain the desired result. For $\eta^-$, we similarly apply \eqref{U1} with $a=\frac{5+2l-m}{4}$ and $b =l+\frac{1}{2}$ to obtain the desired result. 

The assertion about the highest and lowest weights follow from the action of $\eta^\pm$ as differential operators; since the weight vectors that are annihilated by each of these are the ones correspondent to the weights $\mp(2l+1)$ respectively. Directly evaluating and observing that $\ _1F_1(a,a,z)=e^z$ and $\ _1F_1(0,b,z)=1$, the given expressions are obtained.
\end{proof}
\end{proposition}
\begin{definition}
Let $H_l=\ker(\Omega'' - 2\lambda)_K$ denote the $K$-finite vectors in $\ker(\Omega'' - 2\lambda)\subset I''(q,r,s)$. For $k\in \mathbb Z^{\geq 0}$ define
\begin{equation}
	H_k=\spn_{\C}\{\Psi_{m,k}:m\equiv 2k+q \mod 4 \text{ for } m\in \mathbb Z\}. 
	\end{equation}
For $q \equiv 1\mod 4$ and $k\in \mathbb Z^{\geq 0}$ define
\begin{equation}	
H_k^+=\spn_{\C}\{\Psi_{m,k} : m \geq 2k+1 \text{ and } m\equiv 2k+1\mod 4  \text{ for } m\in \mathbb Z\}.
\end{equation}
For $q\equiv -1 \mod 4$ and $k\in \mathbb Z^{\geq 0}$ define
\begin{equation}	
H_k^-=\spn\{\Psi_{m,k} : m \leq -(2k+1) \text{ and } m\equiv -(2k+1)\mod 4 \text{ for } m\in \mathbb Z\}.
\end{equation}
\end{definition}

\begin{lemma}\label{irredLemma}
If $q \equiv 1\mod 4$, then $H_l^+$ is the unique irreducible $\mathfrak{sl}_2$-submodule of $H_l$. If $q \equiv -1\mod 4$, then $H_l^-$ is the unique irreducible $\mathfrak{sl}_2$-submodule of $H_l$.
\begin{proof}
From Equation \eqref{etasAction}, the representation is irreducible whenever $\pm(2l+1)\neq m$ for any $m\equiv 2l+q \mod 4$, this occurs when $q\in 2 \mathbb Z$. We can have $2l+1= m$ for some $m \equiv 2l+q \mod 4$ iff $q \equiv 1 \mod 4$ and a lowest weight occurs in this case. Similarly, $q \equiv -1 \mod 4$ implies that a highest weight occurs. Since the highest and lowest weight cannot occur in the same representation, the action Equation \eqref{etasAction} implies that $H_l^-$ (resp. $H_l^+$) is clearly the unique irreducible submodule of $H_l$. 
\end{proof} 
\end{lemma}

\begin{theorem}
Given $q \in \mathbb Z_4$ and $l=\frac{1}{2}(1+\sqrt{1+8\lambda})$, then as $\mathfrak{sl}_2$-modules:
\begin{enumerate}
 	\item \label{Th1} If $q \equiv 0 \mod 4$ or $q \equiv 2\mod 4$  then $H_l=\ker(\Omega'' - 2\lambda)_K$ is irreducible as an $\mathfrak{sl}_2$-module. 
 	\item \label{Th2} If $q \equiv -1\mod 4$, then $H_l^-$is he only irreducible submodule and the composition series for $\ker(\Omega'' - 2\lambda)_K$ is given by $$0\subset H_l^-\subset H_l$$
 	\item \label{Th3} If $q \equiv 1\mod 4$, then $H_l^+$ is the only irreducible submodule and the composition series of $\ker(\Omega'' - 2\lambda)_K$ is given by $\mathfrak{sl}_2$-submodule $$0\subset H_l^+\subset H_l$$
\end{enumerate}
\begin{proof}
Follows from Lemma \ref{irredLemma}.
\end{proof}
\end{theorem}

\section{Heisenberg action and connections with other kernels}

In this section we will examine the action of the Heisenberg algebra. This will allow us to join all the (non-zero) eigenspaces together in one representation.

Recall that the element $(u,v,0) \in \mathfrak h_3(\R)$ acts on $I'(q,r,s)$ by $(tv-u)\partial_x-2svx$ so, under the isomorphism \eqref{CompactPicIsomorphism}, the elements $$E_\mp:=(1,\pm i, 0 )\in \mathfrak h_3(\C)$$ act by the differential operators
$$\mp e^{\pm i \theta}(\pm \partial_y-2isy)$$
\begin{proposition}\label{ActionsOfHeisOnKTypes}Let $m\in \mathbb Z$ and $k\in \mathbb Z^{\geq 0}$. Then,
$$E^-.\Psi_{m,k}=\frac{(1+2k-m)(k-1)}{(2k-1)(2k+1)}\Psi_{m-2,k+1}-k\Psi_{m-2,k-1}$$
and
$$E^+.\Psi_{m,k}=\frac{(1+2k+m)(k-1)}{2(2k-1)}\Psi_{m+2,k+1}-k\Psi_{m+2,k-1}.$$
\begin{proof}
Combining \eqref{U1} and \eqref{U4} with $a+1$ instead of $a$, one obtains
\begin{equation}\label{Uno}
\,_1F_1(a,b,z)=\,_1F_1(a,b-1,z)-\frac{az}{b(b-1)}\,_1F_1(a+1,b+1,z).
\end{equation}
Using Equation \eqref{U4} with $b+1$ in place of $b$ and combining it with \eqref{U2}, one obtains
\begin{equation}\label{Dos}
\,_1F_1(a,b,z)=-\,_1F_1(a-1,b-1,z)+\frac{b-a}{b-1}z\,_1F_1(a,b+1,z).
\end{equation}
Using Equation \eqref{U0} we can compute the action of $E^\pm$ directly. Let $a = \frac{1+2k-m}{4}$ and $b=k+1/2$. Then it is straightforward to see
\begin{multline*}
E^-.\Psi_{m,k}=-\frac{1}{2}e^{-i(m+2)\theta/2-y^2/2}y^{k-1}\bigl( (-1 + 2 b)\,_1F_1(a,b,y^2)\\ +4a y^2/b\,_1F_1(a+1,b+1,y^2)\bigr). \end{multline*}
Applying \eqref{Uno}, one gets the first equation.

A similar calculation using Equation \eqref{U0} shows
\begin{multline*}
E^+.\Psi_{m,k}=\frac{1}{2}e^{-i(m-2)\theta/2-y^2/2}y^{k-1}\bigl( (1 - 2 b+4y^2)\,_1F_1(a,b,y^2)\\ -4a y^2/b\,_1F_1(a+1,b+1,y^2)\bigr). \end{multline*}
An application of \eqref{U2} gives 
\begin{multline*}
E^+.\Psi_{m,k}=-\frac{1}{2}e^{-i(m-2)\theta/2-y^2/2}y^{k-1}\bigl( 4(b-1)\,_1F_1(a-1,b-1,y^2)\\ +(3-2b)\,_1F_1(a,b,y^2)\bigr) \end{multline*}
and substituting in the expression in \eqref{Dos} gives the desired result.
\end{proof}
\end{proposition}

From Equation \eqref{FormOfK}, it follows that if the eigenvalue $\lambda$ corresponds to the parameter $l$, then $\lambda+l+1$ corresponds to the parameter $l+1$ and $\lambda -l$ corresponds to the parameter $l-1$.

The following corollary is obtained immediately from the previous proposition. It will be useful in seeing that the action of the Heisenberg algebra preserves the structure of highest and lowest weight submodules in $H$.

\begin{corollary}\label{HeisActionOnHiLoWeights}
If $m= 2k+1$, the action of $E^\pm$ on the lowest weight is given by
$$E^-.\psi_{2k+1,k}=-k\Psi_{2k-1,k-1}$$
and
$$E^+.\psi_{2k+1,k}=\frac{(2k+1)(k-1)}{2k-1}\Psi_{2k+3,k+1}-k\Psi_{2k+3,k-1}.$$
If $m=-(2k+1)$, the action of $E^\pm$ on the lowest weight is given by
$$E^-.\psi_{-(2k+1),k}=\frac{2(k-1)}{2k-1}\Psi_{-2k-3,k+1}-k\Psi_{-2k-3,k-1}$$
and
$$E^+.\psi_{-(2k+1),k}=-k\Psi_{-2k+1,k-1}.$$
\begin{proof}
This follows directly from the previous proposition.
\end{proof}
\end{corollary}

We now will show how Proposition \ref{ActionsOfHeisOnKTypes} and Corollary \ref{HeisActionOnHiLoWeights} imply that the action of $\mathfrak h_3$ ties together the kernels indexed by $k$, in a $\mathfrak g$-module. Recall, the cases where $k=0$ and $k=1$ correspond to the potential free case. 
\begin{definition}
Let
\begin{equation}
H =\bigoplus_{l\in \mathbb Z ^{\geq 0}}H_l.
\end{equation}
Whenever the spaces are defined, let
\begin{equation}
H^\pm =\bigoplus_{l\in \mathbb Z^{\geq 2}}H^\pm_l.
\end{equation}
\end{definition}

\begin{theorem}\label{BigTheorem}
Let $q\in \mathbb Z_4$ and $k \in \mathbb Z^{\geq 0}$. With respect to the action of $\mathfrak g$:
\begin{enumerate}
 \item If $q = 0$ or $q = 2$, the composition series of $H$ is $$0\subset H_0\oplus H_1 \subset H.$$
 \item \label{Th5} If $q \equiv -1\mod 4$, then the composition series of $\mathfrak{g}$-submodules of $H$ is as follows  $$0\subset H_0^-\oplus H_1^- \subset H_0\oplus H_1\subset H_0\oplus H_1\oplus H^-\subset H.$$ 
 \item \label{Th6} If $q \equiv 1\mod 4$, then then the composition series of $\mathfrak{g}$-submodules of $H$ is  $$0\subset H_0^+\oplus H_1^+ \subset H_0\oplus H_1\subset H_0\oplus H_1\oplus H^+\subset H.$$ 
 
 \end{enumerate}
\begin{proof}
Let $q \equiv 0 \mod 4$ or $q \equiv 2\mod 4$. Proposition \ref{ActionsOfHeisOnKTypes} shows that the action of $E^\pm$ sends elements in $H_0$ only to $H_1$ and the action of $E^\pm$ sends elements in $H_1$ only to $H_0$. Under this assumption on $q$, each $H_k$ is irreducible under the $\mathfrak{sl}_2$ action, thus $H_0\oplus H_1$ is irreducible under the $\mathfrak g$ action. Now we look at the quotient $H/(H_0\oplus H_1)$. Let $\pi:H \to H/(H_0\oplus H_1)$ be the natural projection. Let $\overline{H}_k$ denote the image of $H_k$ under $\pi$, then the image of $H$ under $\pi$ can be decomposed as a direct sum as $\overline{H} = \bigoplus_j\overline H_{k_j}$ as an $\mathfrak{sl}_2$-module. Proposition \ref{ActionsOfHeisOnKTypes} implies that $E^\pm.\overline H_{k_j}$ has a non-zero component in  $\overline H_{k_j-1}$ and in $\overline H_{k_j+1}$, for $k_j\geq 2$. Since the $\overline H_{k_j-1}$ and $\overline H_{k_j+1}$ are inequivalent $\mathfrak{sl}_2$-representations, then $E^\pm.\overline H_{k_j}$ generates $\overline H_{k_j-1}\oplus \overline H_{k_j+1}$ under the action of $\mathfrak{sl}(2,\R)$. Irreducibility follows easily from this.

Since, the proofs of \eqref{Th5} and \eqref{Th6} are essentially identical, we will only look at the proof of \eqref{Th6}. Irreducibility of $H^+_0\oplus H^+_1$  under $\mathfrak g$ follows from irreducibility under $\mathfrak{sl}_2$ and from the actions on the lowest weights described in Corollary \ref{HeisActionOnHiLoWeights}. 

Next we look at the quotient $(H_0\oplus H_1)/(H^+_0\oplus H^+_1)$. For $j\in \{0,1\}$, write $\overline H_j$ for the image of $H_j$ under the natural projection $H_0\oplus H_1 \to(H_0\oplus H_1)/(H^+_0\oplus H^+_1)$. Then, $\overline H_0$ gets sent to $\overline H_1$ and $\overline H_1$ gets sent to $\overline H_0$ by the action of the Heisenberg algebra. This, together with irreducibility of $\overline H_0$ and $\overline H_1$ under $\mathfrak{sl}_2$, gives irreducibility under $\mathfrak g$.

Finally, we look at the quotient $(H_0\oplus H_1\oplus H^+)/(H_0\oplus H_1)$.  Write $\bigoplus_{k\geq 0} \overline {H}^+_k$ for the image of $H_0\oplus H_1\oplus H^+$ under the natural projection.  The Heisenberg algebra acts as before, and $E^\pm.\overline H^+_{k_j}$ has a  component in  $\overline H^+_{k_j-1}$ and in $\overline H^+_{k_j+1}$, for $k_j \geq 2$. Hence any non-zero element in $\overline H_k$ for $k\geq 2$ generates the whole space. 
\end{proof}
\end{theorem}

Part \eqref{Th5} of Theorem \ref{BigTheorem} can be seen pictorially as follows:

 \begin{center}
\scalebox{1} 
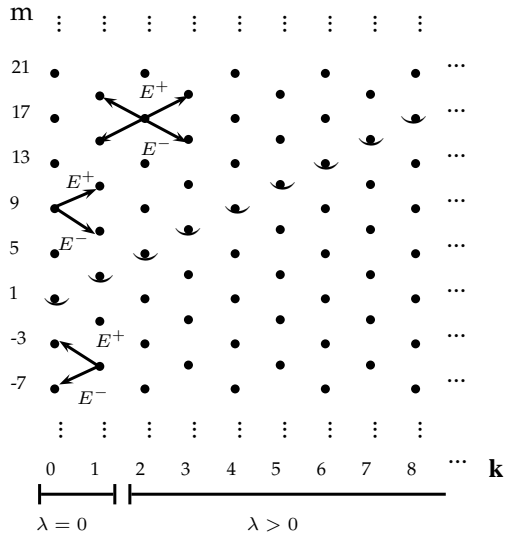
{
\begin{pspicture}(0,-3.5414062)(6.843125,3.5414062)
\usefont{T1}{ppl}{m}{n}
\rput(0.36046875,3.4098437){m}
\psdots[dotsize=0.12](1.396875,-0.68015623)
\psdots[dotsize=0.12](1.396875,-0.08015625)
\psdots[dotsize=0.12](1.396875,0.51984376)
\psdots[dotsize=0.12](1.396875,1.1198437)
\psdots[dotsize=0.12](1.396875,1.7198437)
\psdots[dotsize=0.12](1.396875,2.3198438)
\psdots[dotsize=0.12](1.996875,-1.5801562)
\psdots[dotsize=0.12](1.996875,-0.98015624)
\psdots[dotsize=0.12](1.996875,-0.38015625)
\psdots[dotsize=0.12](1.996875,0.21984375)
\psdots[dotsize=0.12](1.996875,0.81984377)
\psdots[dotsize=0.12](1.996875,1.4198438)
\psdots[dotsize=0.12](1.996875,2.0198438)
\psdots[dotsize=0.12](1.996875,2.6198437)
\psdots[dotsize=0.12](2.576875,-1.2601563)
\psdots[dotsize=0.12](2.576875,-0.66015625)
\psdots[dotsize=0.12](2.576875,-0.06015625)
\psdots[dotsize=0.12](2.576875,0.53984374)
\psdots[dotsize=0.12](2.576875,1.1398437)
\psdots[dotsize=0.12](2.576875,1.7398437)
\psdots[dotsize=0.12](2.576875,2.3398438)
\psdots[dotsize=0.12](3.196875,-1.5801562)
\psdots[dotsize=0.12](3.196875,-0.98015624)
\psdots[dotsize=0.12](3.196875,-0.38015625)
\psdots[dotsize=0.12](3.196875,0.21984375)
\psdots[dotsize=0.12](3.196875,0.81984377)
\psdots[dotsize=0.12](3.196875,1.4198438)
\psdots[dotsize=0.12](3.196875,2.0198438)
\psdots[dotsize=0.12](3.196875,2.6198437)
\psdots[dotsize=0.12](3.796875,-1.2601563)
\psdots[dotsize=0.12](3.796875,-0.66015625)
\psdots[dotsize=0.12](3.796875,-0.06015625)
\psdots[dotsize=0.12](3.796875,0.53984374)
\psdots[dotsize=0.12](3.796875,1.1398437)
\psdots[dotsize=0.12](3.796875,1.7398437)
\psdots[dotsize=0.12](3.796875,2.3398438)
\psdots[dotsize=0.12](4.396875,-1.5801562)
\psdots[dotsize=0.12](4.396875,-0.98015624)
\psdots[dotsize=0.12](4.396875,-0.38015625)
\psdots[dotsize=0.12](4.396875,0.21984375)
\psdots[dotsize=0.12](4.396875,0.81984377)
\psdots[dotsize=0.12](4.396875,1.4198438)
\psdots[dotsize=0.12](4.396875,2.0198438)
\psdots[dotsize=0.12](4.396875,2.6198437)
\psdots[dotsize=0.12](4.996875,-1.2601563)
\psdots[dotsize=0.12](4.996875,-0.66015625)
\psdots[dotsize=0.12](4.996875,-0.06015625)
\psdots[dotsize=0.12](4.996875,0.53984374)
\psdots[dotsize=0.12](4.996875,1.1398437)
\psdots[dotsize=0.12](4.996875,1.7398437)
\psdots[dotsize=0.12](4.996875,2.3398438)
\psdots[dotsize=0.12](5.596875,-1.5801562)
\psdots[dotsize=0.12](5.596875,-0.98015624)
\psdots[dotsize=0.12](5.596875,-0.38015625)
\psdots[dotsize=0.12](5.596875,0.21984375)
\psdots[dotsize=0.12](5.596875,0.81984377)
\psdots[dotsize=0.12](5.596875,1.4198438)
\psdots[dotsize=0.12](5.596875,2.0198438)
\psdots[dotsize=0.12](5.596875,2.6198437)
\usefont{T1}{ppl}{m}{n}
\rput(6.135469,2.7298439){...}
\usefont{T1}{ppl}{m}{n}
\rput(6.135469,2.1298437){...}
\usefont{T1}{ppl}{m}{n}
\rput(6.135469,1.5298438){...}
\usefont{T1}{ppl}{m}{n}
\rput(6.135469,0.9298437){...}
\usefont{T1}{ppl}{m}{n}
\rput(6.135469,0.32984376){...}
\usefont{T1}{ppl}{m}{n}
\rput(6.135469,-0.27015626){...}
\usefont{T1}{ppl}{m}{n}
\rput(6.135469,-0.8701562){...}
\usefont{T1}{ppl}{m}{n}
\rput(6.135469,-1.4701562){...}
\usefont{T1}{ppl}{m}{n}
\rput{-90.0}(-2.4682813,4.0160937){\rput(0.73546875,3.3298438){...}}
\usefont{T1}{ppl}{m}{n}
\rput{-90.0}(-1.8682812,4.6160936){\rput(1.3354688,3.3298438){...}}
\usefont{T1}{ppl}{m}{n}
\rput{-90.0}(-1.2682812,5.216094){\rput(1.9354688,3.3298438){...}}
\usefont{T1}{ppl}{m}{n}
\rput{-90.0}(-0.6682813,5.8160934){\rput(2.5354688,3.3298438){...}}
\usefont{T1}{ppl}{m}{n}
\rput{-90.0}(-0.06828117,6.416094){\rput(3.1354687,3.3298438){...}}
\usefont{T1}{ppl}{m}{n}
\rput{-90.0}(0.53171873,7.0160937){\rput(3.7354689,3.3298438){...}}
\usefont{T1}{ppl}{m}{n}
\rput{-90.0}(1.1317186,7.6160936){\rput(4.335469,3.3298438){...}}
\usefont{T1}{ppl}{m}{n}
\rput{-90.0}(1.7317185,8.216093){\rput(4.9354687,3.3298438){...}}
\usefont{T1}{ppl}{m}{n}
\rput{-90.0}(2.331719,8.816093){\rput(5.5354686,3.3298438){...}}
\usefont{T1}{ppl}{m}{n}
\rput{-90.0}(3.5317187,-0.78390634){\rput(1.3354688,-2.0701563){...}}
\usefont{T1}{ppl}{m}{n}
\rput{-90.0}(4.1317186,-0.18390632){\rput(1.9354688,-2.0701563){...}}
\usefont{T1}{ppl}{m}{n}
\rput{-90.0}(4.731719,0.4160936){\rput(2.5354688,-2.0701563){...}}
\usefont{T1}{ppl}{m}{n}
\rput{-90.0}(5.331719,1.0160937){\rput(3.1354687,-2.0701563){...}}
\usefont{T1}{ppl}{m}{n}
\rput{-90.0}(5.931719,1.6160936){\rput(3.7354689,-2.0701563){...}}
\usefont{T1}{ppl}{m}{n}
\rput{-90.0}(6.5317187,2.2160935){\rput(4.335469,-2.0701563){...}}
\usefont{T1}{ppl}{m}{n}
\rput{-90.0}(7.1317186,2.8160934){\rput(4.9354687,-2.0701563){...}}
\usefont{T1}{ppl}{m}{n}
\rput{-90.0}(7.731719,3.4160938){\rput(5.5354686,-2.0701563){...}}
\usefont{T1}{ppl}{m}{n}
\rput(1.9415625,-2.6351562){\scriptsize 2} \usefont{T1}{ppl}{m}{n} \rput(2.5440626,-2.6351562){\scriptsize 3} \usefont{T1}{ppl}{m}{n} \rput(3.15125,-2.6351562){\scriptsize 4} \usefont{T1}{ppl}{m}{n} \rput(3.7429688,-2.6351562){\scriptsize 5} \usefont{T1}{ppl}{m}{n} \rput(4.347344,-2.6351562){\scriptsize 6} \usefont{T1}{ppl}{m}{n} \rput(4.9425,-2.6351562){\scriptsize 7} \usefont{T1}{ppl}{m}{n} \rput(5.546875,-2.6351562){\scriptsize 8} \usefont{T1}{ppl}{m}{n} \rput(6.155469,-2.5501564){...}
\psdots[dotsize=0.12](0.796875,-1.5801562)
\psdots[dotsize=0.12](0.796875,-0.98015624)
\psdots[dotsize=0.12](0.796875,-0.38015625)
\psdots[dotsize=0.12](0.796875,0.21984375)
\psdots[dotsize=0.12](0.796875,0.81984377)
\psdots[dotsize=0.12](0.796875,1.4198438)
\psdots[dotsize=0.12](0.796875,2.0198438)
\psdots[dotsize=0.12](0.796875,2.6198437)
\usefont{T1}{ppl}{m}{n}
\rput{-90.0}(2.9517188,-1.3639064){\rput(0.7554687,-2.0701563){...}}
\usefont{T1}{ppl}{m}{n}
\rput(0.7475,-2.6351562){\scriptsize 0}
\usefont{T1}{ppl}{m}{n}
\rput(1.3376563,-2.6351562){\scriptsize 1} \usefont{T1}{ppl}{m}{n} \rput(0.8928125,-3.3751562){\scriptsize $\lambda=0$} \psline[linewidth=0.04cm,arrowsize=0.05291667cm 2.0,arrowlength=1.4,arrowinset=0.4]{->}(0.796875,0.83984375)(1.356875,1.0798438)
\usefont{T1}{ppl}{m}{n}
\rput(1.1428125,1.2048438){\scriptsize $E^+$} \psline[linewidth=0.04cm,arrowsize=0.05291667cm 2.0,arrowlength=1.4,arrowinset=0.4]{->}(0.816875,0.81984377)(1.336875,0.47984374)
\usefont{T1}{ppl}{m}{n}
\rput(1.0528125,0.38484374){\scriptsize $E^-$} \psline[linewidth=0.04cm,arrowsize=0.05291667cm 2.0,arrowlength=1.4,arrowinset=0.4]{<->}(1.356875,1.6998438)(2.576875,2.3198438)
\psline[linewidth=0.04cm,arrowsize=0.05291667cm 2.0,arrowlength=1.4,arrowinset=0.4]{<->}(1.436875,2.2998438)(2.556875,1.6998438)
\usefont{T1}{ppl}{m}{n}
\rput(2.1228125,2.4048438){\scriptsize $E^+$} \usefont{T1}{ppl}{m}{n} \rput(2.1528125,1.6448437){\scriptsize $E^-$} \usefont{T1}{ppl}{m}{n} \rput(3.6928124,-3.3751562){\scriptsize $\lambda>0$} \psline[linewidth=0.04cm,tbarsize=0.07055555cm 5.0]{|*-}(1.796875,-2.9801562)(5.996875,-2.9801562)
\psline[linewidth=0.04cm,tbarsize=0.07055555cm 5.0]{|*-|*}(0.596875,-2.9801562)(1.596875,-2.9801562)
\usefont{T1}{ppl}{m}{n}
\rput(0.26296875,0.30484375){\scriptsize 5} \usefont{T1}{ppl}{m}{n} \rput(0.26734376,0.90484375){\scriptsize 9} \usefont{T1}{ppl}{m}{n} \rput(0.349375,1.5048437){\scriptsize 13} \usefont{T1}{ppl}{m}{n} \rput(0.25765625,-0.29515624){\scriptsize 1} \usefont{T1}{ppl}{m}{n} \rput(0.3178125,-0.89515626){\scriptsize -3} \usefont{T1}{ppl}{m}{n} \rput(0.3171875,-1.4951563){\scriptsize -7} \usefont{T1}{ppl}{m}{n} \rput(0.35265625,2.7048438){\scriptsize 21} \usefont{T1}{ppl}{m}{n} \rput(0.34875,2.1048439){\scriptsize 17}
\psdots[dotsize=0.12](1.396875,-1.2801563)
\psline[linewidth=0.04cm,arrowsize=0.05291667cm 2.0,arrowlength=1.4,arrowinset=0.4]{->}(1.356875,-1.2601563)(0.856875,-0.9401562)
\usefont{T1}{ppl}{m}{n}
\rput(1.5428125,-0.91515625){\scriptsize $E^+$} \psline[linewidth=0.04cm,arrowsize=0.05291667cm 2.0,arrowlength=1.4,arrowinset=0.4]{->}(1.416875,-1.2601563)(0.856875,-1.5201563)
\usefont{T1}{ppl}{m}{n}
\rput(1.3128124,-1.6751562){\scriptsize $E^-$} \usefont{T1}{ptm}{m}{n} \rput{89.769325}(0.43840414,-1.2048602){\rput(0.78625,-0.38515624){\large (}} \usefont{T1}{ptm}{m}{n} \rput{89.769325}(1.3160665,-1.4860634){\rput(1.36625,-0.08515625){\large (}} \usefont{T1}{ptm}{m}{n} \rput{89.769325}(2.2136483,-1.7872663){\rput(1.96625,0.21484375){\large (}} \usefont{T1}{ptm}{m}{n} \rput{89.769325}(3.0913913,-2.0285501){\rput(2.52625,0.53484374){\large (}} \usefont{T1}{ptm}{m}{n} \rput{89.769325}(4.048732,-2.3897526){\rput(3.18625,0.83484375){\large (}} \usefont{T1}{ptm}{m}{n} \rput{89.769325}(4.9663134,-2.6710362){\rput(3.78625,1.1548438){\large (}} \usefont{T1}{ptm}{m}{n} \rput{89.769325}(5.803976,-2.9920778){\rput(4.36625,1.4148438){\large (}} \usefont{T1}{ptm}{m}{n} \rput{89.769325}(6.721558,-3.273362){\rput(4.96625,1.7348437){\large (}} \usefont{T1}{ptm}{m}{n} \rput{89.769325}(7.579301,-3.534565){\rput(5.52625,2.0348437){\large (}} \usefont{T1}{ptm}{b}{n} \rput(6.6648436,-2.6301563){k} \end{pspicture} }

 \end{center}
Here, the $\mathfrak{sl}_2$-triple leaves each vertical string invariant, acting by rising and lowering operators and by multiplication by a constant. The operator $E^+$ moves from the $K$-type corresponding to $m$ and $k$ to a linear combination of elements of the $K$-types corresponding to $m+2$ and $k\pm 1$, for $k\geq 2$. Similarly $E^-$ sends $m$ to $m-2$ and $k$ to $k\pm 1$. However, for $k=0$ and $k=1$, the action leaves invariant the direct sum $H_0\oplus H_1$ (c.f. Proposition \ref{ActionsOfHeisOnKTypes}). Each vertical strip has a lowest weight, which is distinguished with an inverted bracket. The action of $E^\pm$ also respects the lowest weight structure as stated in Corollary \ref{HeisActionOnHiLoWeights}.

\bibliographystyle{plain}
\bibliography{mybib}
\end{document}